\documentclass[a4paper]{article}
\usepackage{geometry}
\usepackage{amsmath,amssymb,bm,ulem,url,amsthm, float}
\usepackage{lscape}  
\usepackage{graphicx, rotating}  
\usepackage{tabularx,array,makecell}
\graphicspath{{images/}}
\usepackage{comment}

\newcommand{\CC}{\mathbb{C}}

\newtheorem{lemma}{Lemma}

\newtheorem{thm}[lemma]{Theorem}

\newtheorem{defn}{Definition}

\DeclareMathOperator{\diag}{diag}

\DeclareMathOperator{\GL}{GL}
\DeclareMathOperator{\tr}{tr}

\title{Vector-Valued Invariants Associated with All Irreducible Representations for a Finite Group}
\author{ A. K. M. Selim Reza\footnote{Corresponding author}, Manabu Oura, Masashi Kosuda}

\begin{document}
\date{\empty} 
\maketitle

\textbf{Abstract.}
We investigate the complex reflection group $\mathfrak{G}$ associated with the octahedral group, identified as the ninth entry in the Shephard-Todd classification.
We determine all irreducible representations of $\mathfrak{G}$ and compute the character table.
Moreover, for each representation we compute the module of vector-valued invariants and relate it to the fundamental invariants of the octahedral group.
Additionally, we derive explicit dimension formulas for the corresponding rings of invariants.\\



\section{Introduction}
\label{Intro}

Let $G$ be a finite subgroup of $\GL(n,\CC)$. Write
\[
\mathbf{x}=(x_1,x_2,\dots,x_n)^{\mathsf T},
\qquad
\CC[\mathbf{x}]=\CC[x_1,x_2,\dots,x_n].
\]
The group $G$ acts on $\CC[\mathbf{x}]$ by
\[
g . f(\mathbf{x}) = f(g \mathbf{x}), \quad g \in G
\]

where $g\mathbf{x}$ denotes the usual product of the $n\times n$ matrix $g$ with the column vector $\mathbf{x}$.
A polynomial $f\in \CC[\mathbf{x}]$ is called $G$-invariant if
\[
g\cdot f=f \quad \text{for all } g\in G,
\]
equivalently, if $f(g\mathbf{x})=f(\mathbf{x})$ for all $g\in G$.
The set of all invariants forms a subring
\[
\CC[\mathbf{x}]^{G} := \{f\in \CC[\mathbf{x}] : g\cdot f=f \ \text{for all } g\in G\},
\]
called the invariant ring of $G$.

Among the eight groups listed in Table~II of the Shephard-Todd classification~\cite{shephard-todd} that are derived from the octahedral group, we focus on group number~$9$.
Accordingly, throughout the paper, $\mathfrak{G}$ denotes the Shephard-Todd group $G_9$, generated by
\[
T = \frac{1}{\sqrt{2}}
\begin{pmatrix}
1 & 1 \\
1 & -1
\end{pmatrix},
\qquad
D =
\begin{pmatrix}
1 & 0 \\
0 & i
\end{pmatrix}.
\]
The group $\mathfrak{G}$ has order $192$. Since $\mathfrak{G}$ acts on $\CC^2$, we write
\[
\mathbf{x}=(x,y)^{\mathsf T},
\qquad
\mathfrak{R}:=\CC[x,y]^{\mathfrak{G}}.
\]
The invariant ring, $\mathfrak{R}$ is a polynomial algebra generated by two algebraically independent homogeneous invariants, of degrees $8$ and $24$.
A convenient choice of generators is
\[
\theta := x^8 + 14x^4y^4 + y^8,
\qquad
\varphi := x^{24} + 759x^{16}y^8 + 2576x^{12}y^{12} + 759x^8y^{16} + y^{24}.
\]

Since $\mathfrak{G}$ is closely related to the octahedral group, we also record the classical octahedral invariants appearing in work of Blichfeldt~\cite{miller}:
\[
\Gamma := -x^{5}y + xy^{5}, \qquad
\theta := x^8 + 14x^4y^4 + y^8, \qquad
\Delta := x^{12} - 33x^{8}y^{4} - 33x^{4}y^{8} + y^{12}.
\]
In these terms,
\[
\varphi = \Delta^2 + 66\Gamma^4.
\]

The polynomials $\theta$ and $\varphi$ also have a natural interpretation in coding theory:
they are the homogeneous weight enumerators of the extended binary Hamming code $H_8$ and the extended binary Golay code $G_{24}$, respectively~\cite{sloane}.
By Gleason's theorem~\cite{gleason}, they are $\mathfrak{G}$-invariant and generate $\mathfrak{R}$.
Consequently, the Hilbert--Poincar\'e series of $\mathfrak{R}$ is
\[
\frac{1}{(1-t^8)(1-t^{24})}.
\]

The construction of irreducible representations of the (Shephard--Todd) group of order $96$ was initiated in~\cite{kosuda-oura}, where the structure of its centralizer ring and the corresponding representations were analyzed.
Subsequently, vector-valued invariants for certain representations were computed in~\cite{arX}.
In the present paper we extend this line of work to $\mathfrak{G}=G_9$ and determine all irreducible representations of $\mathfrak{G}$, together with the associated vector-valued invariants.

The paper is organized as follows.
In Section~\ref{IR}, we determine all irreducible representations of $\mathfrak{G}$.
Section~\ref{VVI} collects preliminaries on vector-valued invariants.
In Section~\ref{R}, we compute the vector-valued invariants for each representation and establish relations between invariants in various degrees, in which the fundamental invariants $\Gamma$ and $\Delta$ play a central role.
All computations were carried out using SageMath~\cite{sagemath}.
\section{Irreducible Representations of $\mathfrak{G}$}
\label{IR}

Since $\mathfrak G$ has $32$ conjugacy classes, it has $32$ irreducible complex
representations up to equivalence. Let $(\rho_i,V_i)$ $(i=1,\dots,32)$ be a complete set
of representatives. As usual,
\[
|\mathfrak G|=\sum_{i=1}^{32}(\dim\rho_i)^2,\qquad |\mathfrak G|=192.
\]
Moreover, the degree of an irreducible representation should divide the size of $\mathfrak{G}$. Since
$|\mathfrak G|=3\times 2^6$; hence the possible degrees are
\[
\dim\rho\in\{1,2,3,4,6,8,12\}.
\]

We first determine the linear characters of $\mathfrak G$ (equivalently, its abelianization),
obtaining exactly $8$ one-dimensional irreducibles. Next we construct $8$ pairwise inequivalent irreducible representations of degree $2$ and
$8$ pairwise inequivalent irreducible representations of degree $3$.
Thus we have already exhibited $8+8+8=24$ irreducibles, so there remain $32-24=8$ further
irreducible representations, say $\sigma_1,\dots,\sigma_8$.
Subtracting the known contributions from the sum-of-squares identity gives
\[
80
=192-\bigl(8\cdot 1^2+8\cdot 2^2+8\cdot 3^2\bigr)
=\sum_{j=1}^{8}(\dim\sigma_j)^2.
\]
Since there are no further one-dimensional irreducibles, we have
$\dim\sigma_j\in\{2,3,4,6,8,12\}$.
In fact, $\dim\sigma_j\neq 12$ and $\dim\sigma_j\neq 8$, because $12^2>80$ and
$8^2+7\cdot 2^2=92>80$; similarly $\dim\sigma_j\neq 6$ since $6^2+7\cdot 2^2=64$ would
force the remaining $44$ to be a sum of seven squares from $\{4,9,16\}$, which is impossible.
Hence $\dim\sigma_j\in\{2,3,4\}$ for all $j$.

Write $a,b,c$ for the numbers of $\sigma_j$ of degrees $2,3,4$, respectively. Then
\[
a+b+c=8,\qquad 4a+9b+16c=80.
\]
Subtracting $4(a+b+c)=32$ from the second equation yields
\[
5b+12c=48.
\]
Reducing modulo $5$ gives $2c\equiv 3\pmod 5$, hence $c\equiv 4\pmod 5$.
Since $0\le c\le 8$, we must have $c=4$, and then $5b=0$, so $b=0$ and $a=4$.
Therefore the remaining irreducibles consist of four of degree $2$ and four of degree $4$.

In summary, $\mathfrak G$ has $8$ one-dimensional, $12$ two-dimensional, $8$
three-dimensional, and $4$ four-dimensional irreducible representations.
We construct representatives of all of them explicitly.
\subsection*{One-dimensional representations}

A one-dimensional representation is a homomorphism
$\rho:\mathfrak{G}\to \CC^\times$; since $\mathfrak{G}$ is finite, the image consists of
roots of unity. Writing
\[
\rho(T)=t,\qquad \rho(D)=d \qquad (t,d\in\CC^\times),
\]
the defining relations $T^2=I$ and $D^4=I$ force
\[
t^2=1 \ \Rightarrow\ t=\pm 1,
\qquad
d^4=1 \ \Rightarrow\ d\in\{1,i,-1,-i\}.
\]
These choices yield $2\cdot 4=8$ distinct one-dimensional representations. We enumerate them as follows:
\begin{table}[h!]
\centering
\renewcommand{\arraystretch}{1.1}
\begin{tabular}{c|cc}
\textbf{i} & $\rho_i(T)$ & $\rho_i(D)$ \\
\hline
1 & $ 1$ & $ 1$ \\
2 & $ 1$ & $-1$ \\
3 & $ 1$ & $ i$ \\
4 & $ 1$ & $-i$ \\
5 & $-1$ & $ 1$ \\
6 & $-1$ & $-1$ \\
7 & $-1$ & $ i$ \\
8 & $-1$ & $-i$
\end{tabular}
\caption{All one-dimensional representations of $\mathfrak{G}$}
\end{table}

Here $\rho_1$ is the trivial representation, $\rho_6$ is the representation sending both generators
to $-1$, and $\rho_7$ coincides with the determinant representation since $\det(T)=-1$ and $\det(D)=i$.
Moreover, the remaining representations can be expressed in terms of $\rho_6$ and $\rho_7$:
\[
\rho_2 = \rho_7^{\otimes 2},\quad
\rho_3 = \rho_6 \otimes \rho_7^{\otimes 3},\quad
\rho_4 = \rho_6 \otimes \rho_7,\quad
\rho_5 = \rho_6 \otimes \rho_7^{\otimes 2},\quad
\rho_8 = \rho_7^{\otimes 3}.
\]

\subsection*{Two-dimensional representations}

Let $(\rho_9,V_9)$ be the natural representation, i.e.\ the inclusion
$\mathfrak{G}\hookrightarrow \GL(2,\CC)$ with $T$ and $D$ as in
Section~\ref{Intro}. This representation is irreducible: since $\rho_9(D)$ has distinct
eigenvalues, any $1$-dimensional $\rho_9(D)$-stable subspace is spanned by an eigenvector,
and one checks that no such line is preserved by $\rho_9(T)$. 

Tensoring $\rho_9$ with the $8$ one-dimensional representations produces $8$ faithful irreducible
representations of dimension $2$:
\[
\rho_{8+k} := \rho_k \otimes \rho_9 \qquad (k=1,\dots,8).
\]
Equivalently, for $\varepsilon\in\{\pm 1\}$ and $\eta\in\{1,i,-1,-i\}$ define
\[
\psi_{\varepsilon,\eta}(T)=\varepsilon\,T,\qquad
\psi_{\varepsilon,\eta}(D)=\eta\,D.
\]
Listing these as $\rho_9,\dots,\rho_{16}$ gives Table~\ref{tab:2d-faithful}.

\begin{table}[H]
\centering
\renewcommand{\arraystretch}{1.15}
\begin{tabular}{c|cccc}
\textbf{i} & $(\varepsilon,\eta)$ & $\rho_i$ & $\rho_i(T)$ & $\rho_i(D)$ \\
\hline
9  & $(+1,\,1)$   & $\psi_{+1,1}$     & $T$  & $\diag(1,\,i)$ \\
10 & $(+1,\,i)$   & $\psi_{+1,i}$     & $T$  & $\diag(i,\,-1)$ \\
11 & $(+1,\,-1)$  & $\psi_{+1,-1}$    & $T$  & $\diag(-1,\,-i)$ \\
12 & $(+1,\,-i)$  & $\psi_{+1,-i}$    & $T$  & $\diag(-i,\,1)$ \\
13 & $(-1,\,1)$   & $\psi_{-1,1}$     & $-T$ & $\diag(1,\,i)$ \\
14 & $(-1,\,i)$   & $\psi_{-1,i}$     & $-T$ & $\diag(i,\,-1)$ \\
15 & $(-1,\,-1)$  & $\psi_{-1,-1}$    & $-T$ & $\diag(-1,\,-i)$ \\
16 & $(-1,\,-i)$  & $\psi_{-1,-i}$    & $-T$ & $\diag(-i,\,1)$
\end{tabular}
\caption{Faithful two-dimensional representations}
\label{tab:2d-faithful}
\end{table}

The remaining four two-dimensional irreducible representations will be constructed below
as constituents of a suitable tensor product.

\subsection*{Three-dimensional representations}

Following \cite{kosuda-oura}, we realize a three-dimensional irreducible representation
inside $\rho_9\otimes \rho_9$. Let $\langle \alpha_1,\alpha_2\rangle$ be a basis of $V_9$.
Then
\[
\langle
\alpha_1\otimes \alpha_1,\,
\alpha_1\otimes \alpha_2,\,
\alpha_2\otimes \alpha_1,\,
\alpha_2\otimes \alpha_2
\rangle
\]
is a basis of $V_9\otimes V_9$, and with respect to this basis
\[
\rho_{9}^{\otimes 2}(T)=\frac{1}{2}
\begin{pmatrix}
1&1&1&1\\
1&-1&1&-1\\
1&1&-1&-1\\
1&-1&-1&1
\end{pmatrix},
\qquad
\rho_{9}^{\otimes 2}(D)=\diag(1,i,i,-1).
\]
Consider the vectors
\[
\alpha_1'=\alpha_1\otimes\alpha_1,\quad
\alpha_2'=\alpha_1\otimes\alpha_2+\alpha_2\otimes\alpha_1,\quad
\alpha_3'=\alpha_2\otimes\alpha_2.
\]
Their span is $D$-invariant, and a direct
check shows it is also $T$-invariant. Hence it affords a three-dimensional representation
$(\rho_{21},V_{21})$. With respect to $(\alpha_1',\alpha_2',\alpha_3')$ one finds
\[
\rho_{21}(T)=\frac{1}{2}
\begin{pmatrix}
1&2&1\\
1&0&-1\\
1&-2&1
\end{pmatrix},
\qquad
\rho_{21}(D)=\diag(1,i,-1).
\]
Since $\rho_{21}(D)$ has three distinct eigenvalues, any nonzero proper $D$-stable subspace
is spanned by a nonempty proper subset of $\{\alpha_1',\alpha_2',\alpha_3'\}$; inspection of
the columns of $\rho_{21}(T)$ shows no such subspace is $T$-stable. Thus $\rho_{21}$ is
irreducible.

Tensoring with the one-dimensional representations produces the remaining three-dimensional irreducibles:
\[
\rho_{21+k}:=\rho_{k+1}\otimes \rho_{21}\qquad (k=1,\dots,7),
\]
i.e.
\[
\rho_{22}=\rho_2\otimes\rho_{21},\ \rho_{23}=\rho_3\otimes\rho_{21},\ \dots,\ 
\rho_{28}=\rho_8\otimes\rho_{21}.
\]

\subsection*{Four-dimensional representations}

Next, we extract a four-dimensional irreducible subrepresentation of
$\rho_9\otimes\rho_{21}$.
With respect to the lexicographically ordered basis
$\{\alpha_i\otimes \alpha_j' \mid i=1,2;\ j=1,2,3\}$ of $V_9\otimes V_{21}$, we have
\[
(\rho_9\otimes\rho_{21})(T)=\frac{1}{2\sqrt{2}}
\begin{pmatrix}
1&2&1&1&2&1\\
1&0&-1&1&0&-1\\
1&-2&1&1&-2&1\\
1&2&1&-1&-2&-1\\
1&0&-1&-1&0&1\\
1&-2&1&-1&2&-1
\end{pmatrix},
\qquad
(\rho_9\otimes\rho_{21})(D)=\diag(1,i,-1,i,-1,-i).
\]
Define
\[
\alpha_1''=\alpha_1\otimes\alpha_1',\quad
\alpha_2''=\alpha_1\otimes\alpha_2'+\alpha_2\otimes\alpha_1',\quad
\alpha_3''=\alpha_1\otimes\alpha_3'+\alpha_2\otimes\alpha_2',\quad
\alpha_4''=\alpha_2\otimes\alpha_3'.
\]
The span of $\{\alpha_1'',\alpha_2'',\alpha_3'',\alpha_4''\}$ is $D$-invariant and one checks
it is also $T$-invariant. This yields a four-dimensional representation $(\rho_{29},V_{29})$
with matrices
\[
\rho_{29}(T)=\frac{1}{2\sqrt{2}}
\begin{pmatrix}
1&3&3&1\\
1&1&-1&-1\\
1&-1&-1&1\\
1&-3&3&-1
\end{pmatrix},
\qquad
\rho_{29}(D)=\diag(1,i,-1,-i).
\]
As in the three-dimensional case, $\rho_{29}(D)$ has distinct eigenvalues, so every
$D$-stable subspace is a direct sum of $D$-eigenspaces. Checking images of basis vectors
under $\rho_{29}(T)$ shows that no nontrivial proper $D$-stable subspace is $T$-stable.
Hence $\rho_{29}$ is irreducible.

Tensoring $\rho_{29}$ with one-dimensional representations produces four distinct four-dimensional
irreducibles (the remaining twists are equivalent; this equivalence was verified in
SageMath):
\[
\rho_{30}:=\rho_2\otimes\rho_{29},\qquad
\rho_{31}:=\rho_3\otimes\rho_{29},\qquad
\rho_{32}:=\rho_4\otimes\rho_{29}.
\]

\subsection*{The remaining two-dimensional irreducibles}

Finally, we determine the remaining irreducible representations inside
$(\rho_{9}\otimes \rho_{29},\, V_{9}\otimes V_{29})$. Let
\[
\bigl\{\alpha_i\otimes \alpha_j'' \,\bigm|\, i=1,2,\; j=1,2,3,4\bigr\}
\]
be a basis of $V_{9}\otimes V_{29}$, ordered lexicographically. With respect to this ordered basis, the representation matrices of $T$ and $D$ are:

\[
\begin{aligned}
(\rho_{9} \otimes \rho_{29})(T) &= \frac{1}{4} 
\begin{pmatrix}
1 & 3 & 3 & 1 & 1 & 3 & 3 & 1 \\
1 & 1 & -1 & -1 & 1 & 1 & -1 & -1 \\
1 & -1 & -1 & 1 & 1 & -1 & -1 & 1 \\
1 & -3 & 3 & -1 & 1 & -3 & 3 & -1 \\
1 & 3 & 3 & 1 & -1 & -3 & -3 & 1 \\
1 & 1 & -1 & -1 & -1 & -1 & 1 & 1 \\
1 & -1 & -1 & 1 & -1 & 1 & 1 & -1 \\
1 & -3 & 3 & -1 & -1 & 3 & -3 & 1
\end{pmatrix}, \\
(\rho_{9} \otimes \rho_{29})(D) &= \operatorname{diag}(1, i, -1, -i, i, i, -1, -i, 1).
\end{aligned}
\]

Set
\[
\alpha_1'''=\alpha_1\otimes\alpha_1''+\alpha_2\otimes\alpha_4'',
\qquad
\alpha_2'''=\alpha_1\otimes\alpha_3''+\alpha_2\otimes\alpha_2''.
\]
Then $\langle \alpha_1''',\alpha_2'''\rangle$ is stable under both $T$ and $D$, hence affords a
two-dimensional representation $(\rho_{17},V_{17})$ with
\[
\rho_{17}(T)=\frac{1}{2}
\begin{pmatrix}
1&3\\
1&-1
\end{pmatrix},
\qquad
\rho_{17}(D)=\diag(1,-1).
\]
Since $\rho_{17}(D)$ has two distinct eigenvalues and neither eigenline is $T$-stable,
$\rho_{17}$ is irreducible. Tensoring $\rho_{17}$ with one-dimensional representations yields the
remaining two-dimensional irreducibles (again, the stated
equivalences among twists were verified in SageMath).
\[
\rho_{18}:=\rho_2\otimes\rho_{17},\qquad
\rho_{19}:=\rho_3\otimes\rho_{17},\qquad
\rho_{20}:=\rho_4\otimes\rho_{17}.
\]
At this point we have constructed $32$ irreducible representations:
\[
\underbrace{\rho_1,\dots,\rho_8}_{1\text{-dim}},\quad
\underbrace{\rho_9,\dots,\rho_{20}}_{2\text{-dim}},\quad
\underbrace{\rho_{21},\dots,\rho_{28}}_{3\text{-dim}},\quad
\underbrace{\rho_{29},\dots,\rho_{32}}_{4\text{-dim}}.
\]
Since the number matches the number of conjugacy classes, this is a complete set of irreducible representations of $\mathfrak{G}$. The corresponding character table is given
in the appendix (with $\zeta=e^{i\pi/4}$).
For the presentation of it, we arrange the conjugacy classes into five natural families, represented by elements of the form
\[
\zeta^{k}I,\qquad \zeta^{k}D^{2},\qquad \zeta^{k}D,\qquad \zeta^{k}T,\qquad \zeta^{k}TD.
\]
Accordingly, we reorder the columns of the character table so that the conjugacy classes appear in the order
\[
\{\zeta^{k}I\}_{k=0}^{7},\quad
\{\zeta^{k}D^{2}\}_{k=0}^{3},\quad
\{\zeta^{k}D\}_{k=0}^{7},\quad
\{\zeta^{k}T\}_{k=0}^{3},\quad
\{\zeta^{k}TD\}_{k=0}^{7}.
\]
For each class we record a representative, together with the order of the representative and the size of the conjugacy class.

\section{Vector-valued invariants}
\label{VVI}

Let \(G\) be a finite subgroup of \(GL(n,\CC)\). Put \(V=\CC^n\) with coordinates
\(\mathbf{x}=(x_1,\dots,x_n)\), and let \(S=\CC[\mathbf{x}]\) be the coordinate ring.
The natural action of \(G\) on \(V\) induces an action on \(S\) by substitution:
\[
\sigma\cdot f(\mathbf{x}) \;:=\; f(\sigma \mathbf{x}),
\qquad f\in S,\ \sigma\in G.
\]
With this convention, one has \((\sigma\tau)\cdot f=\tau\cdot(\sigma\cdot f)\).

\medskip

Fix a complex representation \(\rho:G\to GL(m,\CC)\). We view a column vector
\(F(\mathbf{x})=(f_1(\mathbf{x}),\dots,f_m(\mathbf{x}))^{\mathsf{T}}\in S^m\) as an \(m\)-tuple of polynomials.
The action of \(G\) on \(S\) extends componentwise to an action on \(S^m\) by
\[
\sigma\cdot F(\mathbf{x}) \;:=\;
\begin{pmatrix}
\sigma\cdot f_1(\mathbf{x})\\ \vdots\\ \sigma\cdot f_m(\mathbf{x})
\end{pmatrix}
=
\begin{pmatrix}
f_1(\sigma\mathbf{x})\\ \vdots\\ f_m(\sigma\mathbf{x})
\end{pmatrix}.
\]

\begin{defn}
A \(\rho\)-covariant (or vector-valued invariant of \(\rho\))
is a vector \(F\in S^m\) such that
\begin{equation}\label{eq:rho-covariant}
\sigma\cdot F \;=\; \rho(\sigma)\,F
\qquad\text{for all }\sigma\in G,
\end{equation}
where \(\rho(\sigma)F\) denotes the usual matrix multiplication in \(\CC^m\)
(\textit{cf.}~\cite{freitag}).
\end{defn}

Equivalently, \eqref{eq:rho-covariant} can be written as
\[
F(\sigma\mathbf{x}) \;=\; \rho(\sigma)\,F(\mathbf{x})
\qquad (\sigma\in G).
\]
The set of all \(\rho\)-covariants is a graded \(S^G\)-module:
\[
M(\rho)
\;:=\;
\bigl\{\,F\in S^m \ \big|\ F(\sigma\mathbf{x})=\rho(\sigma)F(\mathbf{x})\ \text{for all }\sigma\in G\,\bigr\}.
\]
If \(M(\rho)_k\) denotes the subspace of homogeneous covariants of total degree \(k\),
then each \(M(\rho)_k\) is finite-dimensional over \(\CC\). The graded dimensions are
encoded by the Hilbert series
\[
\Phi_{M(\rho)}(t) \;:=\; \sum_{k\ge 0} \dim_{\CC}\!\bigl(M(\rho)_k\bigr)\,t^k.
\]

\medskip

The Hilbert series is given by the equivariant Molien formula
(see~\cite{gatermann-1996,worfolk-1994}):
\begin{equation}\label{eq:equivariant-molien}
\Phi_{M(\rho)}(t)
\;=\;
\frac{1}{|G|}
\sum_{\sigma\in G}
\frac{\tr\bigl(\rho(\sigma^{-1})\bigr)}{\det\!\bigl(I - t\,\sigma\bigr)}.
\end{equation}
Here \(\det(I-t\,\sigma)\) is computed in the natural \(n\)-dimensional representation of
\(G\) on \(V\), while the trace in the numerator is taken in the representation \(\rho\).

\medskip

In our setting \(n=2\), so the denominator in \eqref{eq:equivariant-molien} is computed from
the natural \(2\times 2\) action of \(\mathfrak{G}\) on \((x,y)\).
Moreover, we have
\(\mathfrak{R}=S^{\mathfrak{G}}=\CC[\theta,\phi]\) with \(\deg\theta=8\) and \(\deg\phi=24\).
Hence every \(\Phi_{M(\rho)}(t)\) is a rational function with denominator
\[
(1-t^{8})(1-t^{24}).
\]

\section{Results}
\label{R}
In this section, we compute all vector-valued invariants explicitly by considering a general function \(F(\mathbf{x})\) with undetermined coefficients. We then impose the condition given in equation~\eqref{eq:rho-covariant} for each representation \(\rho_i \, (i=1,2,\dots 32)\) discussed in Section~2. Solving the resulting system yields a free \(\mathfrak{R}\)-module contained in \(M(\rho_i)\). From these explicit descriptions of the vector-valued invariants, we derive the relations for each representation, which are stated in the following theorems.
 
\begin{thm}
The eight one-dimensional representations are indexed by pairs $(a,b)$ with
$a\in\{0,1,2,3\}$ and $b\in\{0,1\}$, and module \(M(\rho_{a,b})\) is a free \(\mathfrak{R}\)-module of rank $1$ generated by $\Gamma^{a}\Delta^{b}$:
\[
M(\rho_{a,b}) = \mathfrak{R} \,\Gamma^{a}\Delta^{b}.
\]
In particular, for each \( \rho_i \) with \( i=1,2,3,\dots,8 \), the generators are
\[
1,\quad \Gamma^2, \quad \Gamma, \quad \Gamma^3, \quad \Delta, \quad \Gamma^2\Delta, \quad \Gamma\Delta, \quad \Gamma^3\Delta\]
with the corresponding dimension formulas of the form
\[
\frac{t^{k}}{(1 - t^{8})(1 - t^{24})},
\qquad k \in \{0,12,6,18,12,24,18,30\}.
\]
\end{thm}

\begin{thm}
Let $\rho_i$ $(9\le i\le 20)$ be a two-dimensional irreducible representation.
For each $i=9,\dots,20$, the graded $\mathfrak R$-module $M(\rho_i)$ is free of rank $2$.
Moreover, there exist homogeneous elements $\mathbf u_i^{(1)},\mathbf u_i^{(2)}\in M(\rho_i)$ such that
\[
M(\rho_i)=\mathfrak R\,\mathbf u_i^{(1)}\oplus \mathfrak R\,\mathbf u_i^{(2)}.
\]
For a choice of homogeneous generators $(\mathbf u_i^{(1)},\mathbf u_i^{(2)})$ as in \cite{reza},
there exist $c_i\in\mathbb C^\times$ and $k_i\in\mathbb Z$ such that
\[
\det\bigl([\mathbf u_i^{(1)}\ \mathbf u_i^{(2)}]\bigr)=c_i\,\Delta\,\Gamma^{k_i}.
\]
The corresponding values of $(c_i,k_i)$ and the degrees $\deg(\mathbf u_i^{(1)})$, $\deg(\mathbf u_i^{(2)})$
are listed in Table~\ref{tab:gens-2d}.
In particular,
\[
\deg(\mathbf u_i^{(1)})+\deg(\mathbf u_i^{(2)})
=\deg(\Delta)+k_i\deg(\Gamma)=12+6k_i.
\]
\begin{table}[H]
\centering
\[
\renewcommand{\arraystretch}{1.15} 
\begin{array}{c|ccccl}
\rho_i & c_i & k_i & \deg u^{(1)}_i & \deg u^{(2)}_i & \text{Dimension formula (series)}\\
\hline\noalign{\vskip 2pt}
 9  &   -1    & 1 &  1  & 17  & t + t^{9} + 2t^{17} + 3t^{25} + 3t^{33} + \dots\\
10  &  9    & 3 &  7  & 23  & t^7 + t^{15} + 2t^{23} + 3t^{31} + 3t^{39} + \dots\\
11  &   -1    & 5 & 13  & 29  & t^{13} + t^{21} + 2t^{29} + 3t^{37} + 3t^{45} + \dots\\
12  & 176   & 3 & 11  & 19  & t^{11} + 2t^{19} + 2t^{27} + 3t^{35} + 4t^{43} + \dots\\
13  &  4    & 1 &  5  & 13  & t^5 + 2t^{13} + 2t^{21} + 3t^{29} + 4t^{37} + \dots\\
14  &  6    & 3 & 11  & 19  & t^{11} + 2t^{19} + 2t^{27} + 3t^{35} + 4t^{43} + \dots\\
15  &   -1    & 5 & 17  & 25  & t^{17} + 2t^{25} + 2t^{33} + 3t^{41} + 4t^{49} + \dots\\
16  & 1773  & 3 &  7  & 23  & t^7 + t^{15} + 2t^{23} + 3t^{31} + 3t^{39} + \dots\\
17  & -411   & 2 &  8  & 16  & t^8 + 2t^{16} + 2t^{24} + 3t^{32} + 4t^{40} + \dots\\
18  & -10    & 4 & 10  & 26  & t^{10} + t^{18} + 2t^{26} + 3t^{34} + 3t^{42} + \dots\\
19  &   1    & 2 &  4  & 20  & t^4 + t^{12} + 2t^{20} + 3t^{28} + 3t^{36} + \dots\\
20  &  -1    & 4 & 14  & 22  & t^{14} + 2t^{22} + 2t^{30} + 3t^{38} + 4t^{46} + \dots
\end{array}
\]
\caption{Constants $(c_i,k_i)$ and Dimension Formula for $\rho_i$ ($9\le i\le 20$).}
\label{tab:gens-2d}
\end{table}
\end{thm}

\begin{thm}
For each $i\in\{21,\dots,28\}$ let $M(\rho_i)\subset\CC[x,y]^3$ be the covariant
module affording the $3$–dimensional irreducible representation $\rho_i$, generated
by the explicit homogeneous triples
\[
\mathbf{v}^{(1)}_{i},\ \mathbf{v}^{(2)}_{i},\ \mathbf{v}^{(3)}_{i}
\]
listed in \cite{reza}. Then:

\begin{enumerate}

\item
Each $M(\rho_i)$ is a free \(\mathfrak{R}\)–module of rank $3$ with homogeneous basis
$\{\mathbf{v}^{(1)}_{i},\mathbf{v}^{(2)}_{i},\mathbf{v}^{(3)}_{i}\}$.
Writing $\deg\mathbf{v}^{(j)}_{i}$ for the common total degree of its components, the
degrees form an arithmetic progression of step $8$:
\[
\deg\mathbf{v}^{(1)}_{i} = d_i,\qquad
\deg\mathbf{v}^{(2)}_{i} = d_i+8,\qquad
\deg\mathbf{v}^{(3)}_{i} = d_i+16,
\]
where
\[
\begin{array}{c|cccccccc}
i & 21 & 22 & 23 & 24 & 25 & 26 & 27 & 28 \\ \hline
d_i &  2 &  6 &  8 &  4 &  6 & 10 & 12 &  8
\end{array}
\]

\item
The graded Hilbert–Poincar\'e series of $M(\rho_i)$ is
\[
H_{M(\rho_i)}(t)
  = \frac{t^{d_i} + t^{d_i+8} + t^{d_i+16}}{(1 - t^{8})(1 - t^{24})}.
\]
Equivalently, for each degree $d\ge d_i$, the homogeneous component
$(M(\rho_i))_d$ is spanned by the vectors
\[
\theta^{\alpha}\varphi^{\beta}\,\mathbf{v}^{(k)}_{i,j},
\qquad
d = \deg\mathbf{v}^{(j)}_{i} + 8\alpha + 24\beta.
\]

\item Let
\[
D_i:=\det[\mathbf v_{i,1}\ \mathbf v_{i,2}\ \mathbf v_{i,3}]
\]
(the determinant of the $3\times 3$ matrix with columns $\mathbf v_{i,1},\mathbf v_{i,2},\mathbf v_{i,3}$).
Then
\[
D_i \;=\; c_i\,\Delta^{e_i}\,\Gamma^{k_i}
\qquad (c_i\in\mathbb{C}^{\times},\ e_i\in\{1,2\},\ k_i\in\mathbb{Z}),
\]
with $(c_i,e_i,k_i)$ given in the Table~\ref{tab:gens-3d}.
In particular,
\[
\deg D_i=\deg(\mathbf v_{i,1})+\deg(\mathbf v_{i,2})+\deg(\mathbf v_{i,3})
=12e_i+6k_i.
\]
\begin{table}[H]
\centering
\[
\renewcommand{\arraystretch}{1.15} 
\begin{array}{c|cccl}
\rho_i & c_i & e_i & k_i & \text{Dimension formula (series)}\\
\hline\noalign{\vskip 2pt}
21 & -2     & 1 & 3 & t^{2} + 2t^{10} + 3t^{18} + 4t^{26} + 5t^{34} + \dots\\
22 & 4      & 1 & 5 & t^{6} + 2t^{14} + 3t^{22} + 4t^{30} + 5t^{38} + \dots\\
23 & -67    & 1 & 6 & t^{8} + 2t^{16} + 3t^{24} + 4t^{32} + 5t^{40} + \dots\\
24 & 1600   & 1 & 4 & t^{4} + 2t^{12} + 3t^{20} + 4t^{28} + 5t^{36} + \dots\\
25 & -19228 & 2 & 3 & t^{6} + 2t^{14} + 3t^{22} + 4t^{30} + 5t^{38} + \dots\\
26 & 2580   & 2 & 5 & t^{10} + 2t^{18} + 3t^{26} + 4t^{34} + 5t^{42} + \dots\\
27 & 80     & 2 & 6 & t^{12} + 2t^{20} + 3t^{28} + 4t^{36} + 5t^{44} + \dots\\
28 & -894   & 2 & 4 & t^{8} + 2t^{16} + 3t^{24} + 4t^{32} + 5t^{40} + \dots\\
\end{array}
\]
\caption{Constants $(c_i,e_i,k_i)$ and Dimension Formula for $\rho_i$ ($21\le i\le 28$).}
\label{tab:gens-3d}
\end{table}

\item
Let $\tau$ denote the involution on $\CC[x,y]$ given by $(x,y)\mapsto(y,x)$, and
extend it componentwise to $\CC[x,y]^3$.
Then the three components of each homogeneous generator $\mathbf{v}^{(j)}_{i}$ are
related by $\tau$ in a uniform way:

\begin{itemize}
  \item For $i \in \{21,22,25,26\}$, every homogeneous generator can be written in the form
  \[
  \mathbf{v}^{(j)}_{i}(x,y)
    = \bigl(f^{(j)}_{i}(x,y),\ g^{(j)}_{i}(x,y),\ f^{(j)}_{i}(y,x)\bigr)^{\mathsf T},
  \]
  where $f^{(j)}_{i},g^{(j)}_{i}\in\CC[x,y]$ are homogeneous and $g^{(j)}_{i}$ is symmetric,
  i.e.\ $g^{(j)}_{i}(y,x) = g^{(j)}_{i}(x,y)$.

  \item For $i \in \{23,24,27,28\}$, every homogeneous generator can be written in the form
  \[
  \mathbf{v}^{(j)}_{i}(x,y)
    = \bigl(f^{(j)}_{i}(x,y),\ g^{(j)}_{i}(x,y),\ -\,f^{(j)}_{i}(y,x)\bigr)^{\mathsf T},
  \]
  where $g^{(j)}_{i}$ is skew–symmetric, i.e.\ $g^{(j)}_{i}(y,x) = -\,g^{(j)}_{i}(x,y)$.
\end{itemize}

\end{enumerate}
\end{thm}

\begin{thm}
For each $i\in\{29,30,31,32\}$ let $M(\rho_i)\subset\CC[x,y]^4$ be the covariant
module affording the $4$–dimensional irreducible representation $\rho_i$, generated
by the explicit homogeneous $4$–tuples
\[
\mathbf{w}^{(1)}_{i},\ \mathbf{w}^{(2)}_{i},\ \mathbf{w}^{(3)}_{i},\ \mathbf{w}^{(4)}_{i}
\]
given in \cite{reza}. Then:
\begin{enumerate}
\item
Each $M(\rho_i)$ is a free \(\mathfrak{R}\)–module of rank $4$ with homogeneous basis
$\{\mathbf{w}^{(1)}_{i},\mathbf{w}^{(2)}_{i},\mathbf{w}^{(3)}_{i},\mathbf{w}^{(4)}_{i}\}$.
Writing $\deg\mathbf{w}^{(k)}_{i}$ for the common total degree of its components, the
degrees and dimension formula are:
\begin{table}[H]
\centering
\[
\renewcommand{\arraystretch}{1.15} 
\begin{array}{c|cl}
\rho_i & \deg(\mathbf{w}^{(1)}_{i},\mathbf{w}^{(2)}_{i},\mathbf{w}^{(3)}_{i},\mathbf{w}^{(4)}_{i}) & \text{Dimension formula (series)}\\
\hline\noalign{\vskip 2pt}

29 & (3,11,19,27) 
   & t^3 + 2t^{11} + 3t^{19} + 5t^{27} + 6t^{35} + \dots \\[6pt]
30 & (7,15,15,23) 
   & t^7 + 3t^{15} + 4t^{23} + 5t^{31} + 7t^{39} + \dots \\[6pt]
31 & (9,9,17,25) 
   & 2t^9 + 3t^{17} + 4t^{25} + 6t^{33} + 7t^{41} + \dots \\[6pt]
32 & (5,13,21,21) 
   &  t^5 + 2t^{13} + 4t^{21} + 5t^{29} + 6t^{37} + \dots
\end{array}
\]
\caption{Dimension Formula for $\rho_i$ ($29\le i\le 32$).}
\label{tab:gens-4d}
\end{table}

\item
Set $J:=\Delta\,\Gamma^3$ (so $\deg J=30$). For each $i$ let
\[
D_i:=\det[\mathbf{w}^{(1)}_{i},\mathbf{w}^{(2)}_{i},\mathbf{w}^{(3)}_{i},\mathbf{w}^{(4)}_{i}],
\]
the determinant of the $4\times 4$ matrix with columns $\mathbf{w}^{(1)}_{i},,\dots,\mathbf{w}^{(4)}_{i}$
in the displayed order. Then
\[
D_i \;=\; c_i\,J^2 \;=\; c_i\,\Delta^2\,\Gamma^6
\]
with
\[
c_{29}=3276,\qquad
c_{30}=-11780,\qquad
c_{31}=344960,\qquad
c_{32}=437400.
\]
In particular,
\[
\deg D_i=\deg(\mathbf{w}^{(1)}_{i})+\deg(\mathbf{w}^{(2)}_{i})+\deg(\mathbf{w}^{(3)}_{i})+\deg(\mathbf{w}^{(4)}_{i})=60.
\]
\item
Let $\tau$ denote the involution on $\CC[x,y]$ given by $(x,y)\mapsto(y,x)$, and
extend it component-wise to $\CC[x,y]^4$.
Then the components of each homogeneous generator $\mathbf{w}^{(k)}_{i}$ are related
by $\tau$ in a uniform way:

\begin{itemize}
  \item For $i=29,30$, every homogeneous generator can be written in the form
  \[
  \mathbf{w}^{(k)}_{i}(x,y)
  = \bigl(f^{(k)}_{i}(x,y),\ g^{(k)}_{i}(x,y),\ g^{(k)}_{i}(y,x),\ f^{(k)}_{i}(y,x)\bigr)^{\mathsf T}
  \]
  for suitable homogeneous polynomials $f^{(k)}_{i},g^{(k)}_{i}\in\CC[x,y]$.

  \item For $i=31,32$, every homogeneous generator can be written in the form
  \[
  \mathbf{w}^{(k)}_{i}(x,y)
  = \bigl(f^{(k)}_{i}(x,y),\ g^{(k)}_{i}(x,y),\ -g^{(k)}_{i}(y,x),\ -f^{(k)}_{i}(y,x)\bigr)^{\mathsf T}.
  \]
\end{itemize}

\end{enumerate}

\end{thm}
\medskip
\textbf{Acknowledgements}.
The first named author was supported in part by funds from  the Ministry of Education, Culture, Sports, Science and Technology (MEXT), Japan, and the second and third named authors were supported by JSPS KAKENHI Grant Numbers 24K06827 and 24K06644 respectively.

Department of Mathematics, Khulna University of Engineering \& Technology, Khulna-9203, Bangladesh and
Graduate School of Natural Science and Technology, Kanazawa University, Ishikawa 920-1192, Japan

\textit{Email address}: \texttt{selim\_1992@math.kuet.ac.bd}

\bigskip

Faculty of Mathematics and Physics,
Institute of Science and Engineering,
Kanazawa University,
Kakuma-machi,
Ishikawa 920-1192,
Japan

\textit{Email address}: \texttt{oura@se.kanazawa-u.ac.jp}

\bigskip

Faculty of Engineering, University of Yamanashi, 400-8511, Japan

\textit{Email address}: \texttt{mkosuda@yamanashi.ac.jp}
\pagebreak

\begin{landscape}
\renewcommand{\arraystretch}{1.15}
\resizebox{1.15\textwidth}{!}{$
\begin{array}{c|*{16}{c}}
  & \mathfrak{C}_{01} & \mathfrak{C}_{02} & \mathfrak{C}_{03} & \mathfrak{C}_{04} & \mathfrak{C}_{05} & \mathfrak{C}_{06} & \mathfrak{C}_{07} & \mathfrak{C}_{08} & \mathfrak{C}_{09} & \mathfrak{C}_{10} & \mathfrak{C}_{11} & \mathfrak{C}_{12} & \mathfrak{C}_{13} & \mathfrak{C}_{14} & \mathfrak{C}_{15} & \mathfrak{C}_{16}\\
\hline
\mathrm{Rep.}
& I & \zeta I & \zeta^{2}I & \zeta^{3}I & \zeta^{4}I & \zeta^{5}I & \zeta^{6}I & \zeta^{7}I
& D^{2} & \zeta D^{2} & \zeta^{2}D^{2} & \zeta^{3}D^{2}
& D & \zeta D & \zeta^{2}D & \zeta^{3}D\\
\mathrm{ord}
& 1 & 8 & 4 & 8 & 2 & 8 & 4 & 8
& 2 & 8 & 4 & 8
& 4 & 8 & 4 & 8\\
|\mathcal C|
& 1 & 1 & 1 & 1 & 1 & 1 & 1 & 1
& 6 & 6 & 6 & 6
& 6 & 6 & 6 & 6\\
\hline
\chi_1& 1& 1& 1& 1& 1& 1& 1& 1& 1& 1& 1& 1& 1& 1& 1& 1\\
\chi_2& 1& -1& 1& -1& 1& -1& 1& -1& 1& -1& 1& -1& -1& 1& -1& 1\\
\chi_3& 1& -\zeta^2& -1& \zeta^2& 1& -\zeta^2& -1& \zeta^2& -1& \zeta^2& 1& -\zeta^2& \zeta^2& 1& -\zeta^2& -1\\
\chi_4& 1& \zeta^2& -1& -\zeta^2& 1& \zeta^2& -1& -\zeta^2& -1& -\zeta^2& 1& \zeta^2& -\zeta^2& 1& \zeta^2& -1\\
\chi_5& 1& -1& 1& -1& 1& -1& 1& -1& 1& -1& 1& -1& 1& -1& 1& -1\\
\chi_6& 1& 1& 1& 1& 1& 1& 1& 1& 1& 1& 1& 1& -1& -1& -1& -1\\
\chi_7& 1& \zeta^2& -1& -\zeta^2& 1& \zeta^2& -1& -\zeta^2& -1& -\zeta^2& 1& \zeta^2& \zeta^2& -1& -\zeta^2& 1\\
\chi_8& 1& -\zeta^2& -1& \zeta^2& 1& -\zeta^2& -1& \zeta^2& -1& \zeta^2& 1& -\zeta^2& -\zeta^2& -1& \zeta^2& 1\\
\chi_{9}& 2& 2\zeta& 2\zeta^2& 2\zeta^3& -2& -2\zeta& -2\zeta^2& -2\zeta^3& 0& 0& 0& 0& \zeta^2+1& \zeta^3+\zeta& \zeta^2-1& \zeta^3-\zeta\\
\chi_{10}& 2& -2\zeta^3& -2\zeta^2& -2\zeta& -2& 2\zeta^3& 2\zeta^2& 2\zeta& 0& 0& 0& 0& \zeta^2-1& \zeta^3+\zeta& \zeta^2+1& -\zeta^3+\zeta\\
\chi_{11}& 2& -2\zeta& 2\zeta^2& -2\zeta^3& -2& 2\zeta& -2\zeta^2& 2\zeta^3& 0& 0& 0& 0& -\zeta^2-1& \zeta^3+\zeta& -\zeta^2+1& \zeta^3-\zeta\\
\chi_{12}& 2& 2\zeta^3& -2\zeta^2& 2\zeta& -2& -2\zeta^3& 2\zeta^2& -2\zeta& 0& 0& 0& 0& -\zeta^2+1& \zeta^3+\zeta& -\zeta^2-1& -\zeta^3+\zeta\\
\chi_{13}& 2& -2\zeta& 2\zeta^2& -2\zeta^3& -2& 2\zeta& -2\zeta^2& 2\zeta^3& 0& 0& 0& 0& \zeta^2+1& -\zeta^3-\zeta& \zeta^2-1& -\zeta^3+\zeta\\
\chi_{14}& 2& 2\zeta^3& -2\zeta^2& 2\zeta& -2& -2\zeta^3& 2\zeta^2& -2\zeta& 0& 0& 0& 0& \zeta^2-1& -\zeta^3-\zeta& \zeta^2+1& \zeta^3-\zeta\\
\chi_{15}& 2& 2\zeta& 2\zeta^2& 2\zeta^3& -2& -2\zeta& -2\zeta^2& -2\zeta^3& 0& 0& 0& 0& -\zeta^2-1& -\zeta^3-\zeta& -\zeta^2+1& -\zeta^3+\zeta\\
\chi_{16}& 2& -2\zeta^3& -2\zeta^2& -2\zeta& -2& 2\zeta^3& 2\zeta^2& 2\zeta& 0& 0& 0& 0& -\zeta^2+1& -\zeta^3-\zeta& -\zeta^2-1& \zeta^3-\zeta\\
\chi_{17}& 2& 2& 2& 2& 2& 2& 2& 2& 2& 2& 2& 2& 0& 0& 0& 0\\
\chi_{18}& 2& 2\zeta^2& -2& -2\zeta^2& 2& 2\zeta^2& -2& -2\zeta^2& -2& -2\zeta^2& 2& 2\zeta^2& 0& 0& 0& 0\\
\chi_{19}& 2& -2& 2& -2& 2& -2& 2& -2& 2& -2& 2& -2& 0& 0& 0& 0\\
\chi_{20}& 2& -2\zeta^2& -2& 2\zeta^2& 2& -2\zeta^2& -2& 2\zeta^2& -2& 2\zeta^2& 2& -2\zeta^2& 0& 0& 0& 0\\
\chi_{21}& 3& 3\zeta^2& -3& -3\zeta^2& 3& 3\zeta^2& -3& -3\zeta^2& 1& \zeta^2& -1& -\zeta^2& \zeta^2& -1& -\zeta^2& 1\\
\chi_{22}& 3& -3\zeta^2& -3& 3\zeta^2& 3& -3\zeta^2& -3& 3\zeta^2& 1& -\zeta^2& -1& \zeta^2& -\zeta^2& -1& \zeta^2& 1\\
\chi_{23}& 3& 3& 3& 3& 3& 3& 3& 3& -1& -1& -1& -1& -1& -1& -1& -1\\
\chi_{24}& 3& -3& 3& -3& 3& -3& 3& -3& -1& 1& -1& 1& 1& -1& 1& -1\\
\chi_{25}& 3& -3\zeta^2& -3& 3\zeta^2& 3& -3\zeta^2& -3& 3\zeta^2& 1& -\zeta^2& -1& \zeta^2& \zeta^2& 1& -\zeta^2& -1\\
\chi_{26}& 3& 3\zeta^2& -3& -3\zeta^2& 3& 3\zeta^2& -3& -3\zeta^2& 1& \zeta^2& -1& -\zeta^2& -\zeta^2& 1& \zeta^2& -1\\
\chi_{27}& 3& -3& 3& -3& 3& -3& 3& -3& -1& 1& -1& 1& -1& 1& -1& 1\\
\chi_{28}& 3& 3& 3& 3& 3& 3& 3& 3& -1& -1& -1& -1& 1& 1& 1& 1\\
\chi_{29}& 4& -4\zeta^3& -4\zeta^2& -4\zeta& -4& 4\zeta^3& 4\zeta^2& 4\zeta& 0& 0& 0& 0& 0& 0& 0& 0\\
\chi_{30}& 4& 4\zeta& 4\zeta^2& 4\zeta^3& -4& -4\zeta& -4\zeta^2& -4\zeta^3& 0& 0& 0& 0& 0& 0& 0& 0\\
\chi_{31}& 4& 4\zeta^3& -4\zeta^2& 4\zeta& -4& -4\zeta^3& 4\zeta^2& -4\zeta& 0& 0& 0& 0& 0& 0& 0& 0\\
\chi_{32}& 4& -4\zeta& 4\zeta^2& -4\zeta^3& -4& 4\zeta& -4\zeta^2& 4\zeta^3& 0& 0& 0& 0& 0& 0& 0& 0
\end{array}
$}
\end{landscape}

\begin{landscape}
\renewcommand{\arraystretch}{1.15}
\resizebox{1.15\textwidth}{!}{$
\begin{array}{c|*{16}{c}}
 & \mathfrak{C}_{17} & \mathfrak{C}_{18} & \mathfrak{C}_{19} & \mathfrak{C}_{20} & \mathfrak{C}_{21} & \mathfrak{C}_{22} & \mathfrak{C}_{23} & \mathfrak{C}_{24} & \mathfrak{C}_{25} & \mathfrak{C}_{26} & \mathfrak{C}_{27} & \mathfrak{C}_{28} & \mathfrak{C}_{29} & \mathfrak{C}_{30} & \mathfrak{C}_{31} & \mathfrak{C}_{32}\\
\hline
\mathrm{Rep.}
& \zeta^{4}D & \zeta^{5}D & \zeta^{6}D & \zeta^{7}D
& T & \zeta T & \zeta^{2}T & \zeta^{3}T
& TD & \zeta TD & \zeta^{2}TD & \zeta^{3}TD & \zeta^{4}TD & \zeta^{5}TD & \zeta^{6}TD & \zeta^{7}TD\\
\mathrm{ord}
& 4 & 8 & 4 & 8
& 2 & 8 & 4 & 8
& 24 & 6 & 24 & 12 & 24 & 3 & 24 & 12\\
|\mathcal C|
& 6 & 6 & 6 & 6
& 12 & 12 & 12 & 12
& 8 & 8 & 8 & 8 & 8 & 8 & 8 & 8\\
\hline
\chi_{1} & 1 & 1 & 1 & 1 & 1 & 1 & 1 & 1 & 1 & 1 & 1 & 1 & 1 & 1 & 1 & 1\\
\chi_{2} & -1 & 1 & -1 & 1 & 1 & -1 & 1 & -1 & -1 & 1 & -1 & 1 & -1 & 1 & -1 & 1\\
\chi_{3} & \zeta^2 & 1 & -\zeta^2 & -1 & 1 & -\zeta^2 & -1 & \zeta^2 & \zeta^2 & 1 & -\zeta^2 & -1 & \zeta^2 & 1 & -\zeta^2 & -1\\
\chi_{4} & -\zeta^2 & 1 & \zeta^2 & -1 & 1 & \zeta^2 & -1 & -\zeta^2 & -\zeta^2 & 1 & \zeta^2 & -1 & -\zeta^2 & 1 & \zeta^2 & -1\\
\chi_{5} & 1 & -1 & 1 & -1 & -1 & 1 & -1 & 1 & -1 & 1 & -1 & 1 & -1 & 1 & -1 & 1\\
\chi_{6} & -1 & -1 & -1 & -1 & -1 & -1 & -1 & -1 & 1 & 1 & 1 & 1 & 1 & 1 & 1 & 1\\
\chi_{7} & \zeta^2 & -1 & -\zeta^2 & 1 & -1 & -\zeta^2 & 1 & \zeta^2 & -\zeta^2 & 1 & \zeta^2 & -1 & -\zeta^2 & 1 & \zeta^2 & -1\\
\chi_{8} & -\zeta^2 & -1 & \zeta^2 & 1 & -1 & \zeta^2 & 1 & -\zeta^2 & \zeta^2 & 1 & -\zeta^2 & -1 & \zeta^2 & 1 & -\zeta^2 & -1\\
\chi_{9} & -\zeta^2-1 & -\zeta^3-\zeta & -\zeta^2+1 & -\zeta^3+\zeta & 0 & 0 & 0 & 0 & -\zeta^3 & 1 & \zeta & \zeta^2 & \zeta^3 & -1 & -\zeta & -\zeta^2\\
\chi_{10} & -\zeta^2+1 & -\zeta^3-\zeta & -\zeta^2-1 & \zeta^3-\zeta & 0 & 0 & 0 & 0 & \zeta & 1 & -\zeta^3 & -\zeta^2 & -\zeta & -1 & \zeta^3 & \zeta^2\\
\chi_{11} & \zeta^2+1 & -\zeta^3-\zeta & \zeta^2-1 & -\zeta^3+\zeta & 0 & 0 & 0 & 0 & \zeta^3 & 1 & -\zeta & \zeta^2 & -\zeta^3 & -1 & \zeta & -\zeta^2\\
\chi_{12} & \zeta^2-1 & -\zeta^3-\zeta & \zeta^2+1 & \zeta^3-\zeta & 0 & 0 & 0 & 0 & -\zeta & 1 & \zeta^3 & -\zeta^2 & \zeta & -1 & -\zeta^3 & \zeta^2\\
\chi_{13} & -\zeta^2-1 & \zeta^3+\zeta & -\zeta^2+1 & \zeta^3-\zeta & 0 & 0 & 0 & 0 & \zeta^3 & 1 & -\zeta & \zeta^2 & -\zeta^3 & -1 & \zeta & -\zeta^2\\
\chi_{14} & -\zeta^2+1 & \zeta^3+\zeta & -\zeta^2-1 & -\zeta^3+\zeta & 0 & 0 & 0 & 0 & -\zeta & 1 & \zeta^3 & -\zeta^2 & \zeta & -1 & -\zeta^3 & \zeta^2\\
\chi_{15} & \zeta^2+1 & \zeta^3+\zeta & \zeta^2-1 & \zeta^3-\zeta & 0 & 0 & 0 & 0 & -\zeta^3 & 1 & \zeta & \zeta^2 & \zeta^3 & -1 & -\zeta & -\zeta^2\\
\chi_{16} & \zeta^2-1 & \zeta^3+\zeta & \zeta^2+1 & -\zeta^3+\zeta & 0 & 0 & 0 & 0 & \zeta & 1 & -\zeta^3 & -\zeta^2 & -\zeta & -1 & \zeta^3 & \zeta^2\\
\chi_{17} & 0 & 0 & 0 & 0 & 0 & 0 & 0 & 0 & -1 & -1 & -1 & -1 & -1 & -1 & -1 & -1\\
\chi_{18} & 0 & 0 & 0 & 0 & 0 & 0 & 0 & 0 & \zeta^2 & -1 & -\zeta^2 & 1 & \zeta^2 & -1 & -\zeta^2 & 1\\
\chi_{19} & 0 & 0 & 0 & 0 & 0 & 0 & 0 & 0 & 1 & -1 & 1 & -1 & 1 & -1 & 1 & -1\\
\chi_{20} & 0 & 0 & 0 & 0 & 0 & 0 & 0 & 0 & -\zeta^2 & -1 & \zeta^2 & 1 & -\zeta^2 & -1 & \zeta^2 & 1\\
\chi_{21} & \zeta^2 & -1 & -\zeta^2 & 1 & 1 & \zeta^2 & -1 & -\zeta^2 & 0 & 0 & 0 & 0 & 0 & 0 & 0 & 0\\
\chi_{22} & -\zeta^2 & -1 & \zeta^2 & 1 & 1 & -\zeta^2 & -1 & \zeta^2 & 0 & 0 & 0 & 0 & 0 & 0 & 0 & 0\\
\chi_{23} & -1 & -1 & -1 & -1 & 1 & 1 & 1 & 1 & 0 & 0 & 0 & 0 & 0 & 0 & 0 & 0\\
\chi_{24} & 1 & -1 & 1 & -1 & 1 & -1 & 1 & -1 & 0 & 0 & 0 & 0 & 0 & 0 & 0 & 0\\
\chi_{25} & \zeta^2 & 1 & -\zeta^2 & -1 & -1 & \zeta^2 & 1 & -\zeta^2 & 0 & 0 & 0 & 0 & 0 & 0 & 0 & 0\\
\chi_{26} & -\zeta^2 & 1 & \zeta^2 & -1 & -1 & -\zeta^2 & 1 & \zeta^2 & 0 & 0 & 0 & 0 & 0 & 0 & 0 & 0\\
\chi_{27} & -1 & 1 & -1 & 1 & -1 & 1 & -1 & 1 & 0 & 0 & 0 & 0 & 0 & 0 & 0 & 0\\
\chi_{28} & 1 & 1 & 1 & 1 & -1 & -1 & -1 & -1 & 0 & 0 & 0 & 0 & 0 & 0 & 0 & 0\\
\chi_{29} & 0 & 0 & 0 & 0 & 0 & 0 & 0 & 0 & -\zeta & -1 & \zeta^3 & \zeta^2 & \zeta & 1 & -\zeta^3 & -\zeta^2\\
\chi_{30} & 0 & 0 & 0 & 0 & 0 & 0 & 0 & 0 & \zeta^3 & -1 & -\zeta & -\zeta^2 & -\zeta^3 & 1 & \zeta & \zeta^2\\
\chi_{31} & 0 & 0 & 0 & 0 & 0 & 0 & 0 & 0 & \zeta & -1 & -\zeta^3 & \zeta^2 & -\zeta & 1 & \zeta^3 & -\zeta^2\\
\chi_{32} & 0 & 0 & 0 & 0 & 0 & 0 & 0 & 0 & -\zeta^3 & -1 & \zeta & -\zeta^2 & \zeta^3 & 1 & -\zeta & \zeta^2
\end{array}
$}
\end{landscape}

\end{document}